\documentclass[letterpaper, 10 pt, conference]{ieeeconf}  % Comment this line out if you need a4paper
\IEEEoverridecommandlockouts
\usepackage[top=57pt,bottom=43pt,left=48pt,right=48pt]{geometry}
\usepackage{cite}
\usepackage{amsmath,amssymb,amsfonts}
\usepackage{algorithm}
\usepackage{booktabs} % For prettier tables
\usepackage{algpseudocode}
\usepackage{graphicx}
\usepackage{textcomp}
\usepackage{xcolor}
\usepackage{subfiles}
\usepackage{mathtools}
\usepackage{subcaption}
\usepackage{graphicx}
\usepackage{svg}
\usepackage{url}
\usepackage{hyperref}
\usepackage{comment}
\usepackage{bm}
\newcommand\norm[1]{\lVert#1\rVert}

\newtheorem{theorem}{Theorem}

\newtheorem{lemma}{Lemma}

\newtheorem{assumption}{Assumption}
\DeclareMathOperator*{\argmin}{arg\,min}
\def\BibTeX{{\rm B\kern-.05em{\sc i\kern-.025em b}\kern-.08em
    T\kern-.1667em\lower.7ex\hbox{E}\kern-.125emX}}

\newcommand{\vishnu}[1]{#1}

\begin{document}
%%% Title
\title{A Dynamically Weighted ADMM Framework for Byzantine Resilience}
\author{Vishnu Vijay,
Kartik A. Pant,
Minhyun Cho,
and Inseok Hwang% <-this % stops a space
\thanks{\textcolor{blue}{This work has been submitted to the IEEE for possible publication. Copyright may be transferred without notice, after which this version may no longer be accessible.}}
\thanks{The authors are with the School of Aeronautics and Astronautics, Purdue University,
West Lafayette, IN 47907. Email:
\{
{\tt\small vvijay},
{\tt\small kpant},
{\tt\small cho515},
{\tt\small ihwang}
\}
{\tt\small @purdue.edu}}%
}

\maketitle

%%% Abstract
\begin{abstract}
    The alternating direction of multipliers method (ADMM) is a popular method to solve distributed consensus optimization utilizing efficient communication among various nodes in the network. However, in the presence of faulty or attacked nodes, even a small perturbation (or sharing false data) during the communication can lead to divergence of the solution.
    To address this issue, in this work we consider ADMM under the effect of Byzantine threat, where an unknown subset of nodes is subject to Byzantine attacks or faults.
    We propose Dynamically Weighted ADMM (DW-ADMM), a novel variant of ADMM that uses dynamic weights on the edges of the network, thus promoting resilient distributed optimization. We establish that the proposed method (i) produces a nearly identical solution to conventional ADMM in the error-free case, and (ii) guarantees a bounded solution with respect to the global minimizer, even under Byzantine threat. Finally, we demonstrate the effectiveness of our proposed algorithm using an illustrative numerical simulation.
\end{abstract}


%%% Introduction
\section{Introduction}
\label{sec:intro}
Prompted by recent developments in inexpensive high-performance computational hardware, a great deal of interest has shifted from centralized networked systems, which are significantly constrained by communication and computational costs, to distributed networked systems. A popular problem with these systems is distributed optimization, which involves a group of nodes, each with a local objective function, that seek to collaboratively minimize (or maximize) an aggregate of the local objective functions by sharing only local update information.
This problem is widely popular since many problems in control theory, signal processing, and even machine learning can be transformed into the general framework of distributed optimization \cite{nedic2018distributed, yang2019survey, rabbat2004distributed}. Much research has been conducted in this area, leading to the development of a number of tools and techniques, such as distributed Newton-Raphson, distributed gradient descent, and the alternating direction of multiplier method (ADMM). ADMM has stood out due to its fast convergence speed \cite{boyd2011distributed, shi2014linear}, while also being well suited for a wide variety of distributed computing scenarios \cite{vijay2025range, shethia2021distributed, he2025straggler, li2022robust, ling2016weighted}.
However, due to its dependence on accurate and consistent communication and computation, it is highly susceptible to faults/attacks, where even a slight error in the shared information or computed solution could cause divergence. This necessitates the development of a resilient and secure distributed ADMM algorithm. 

In this work, we explore the situation in which an unknown subset of network nodes suffers Byzantine faults/attacks. This threat model is highly effective against typical distributed optimization algorithms, where even a single affected node could disrupt the network. The Byzantine threat model is similar to a traditional adversarial threat model for networks, but differs in that an affected node can send different messages to each of its neighbors \cite{su2020byzantine}. This difference is marked in the development of consensus algorithms resilient to these threat models. It was shown in \cite{sundaram2018distributed} that a consensus algorithm can tolerate up to half of the network nodes acting adversarially, but cannot tolerate more than $\frac{1}{3}$ of the network nodes being Byzantine. Byzantine threats are not isolated to only distributed optimization, and have been considered for many distributed systems applications. In \cite{zegers2021event}, researchers explored formation control and leader tracking under Byzantine threat, using a reputation system to identify the Byzantine agents and selectively ignore their information. The authors in \cite{cui2024resilient} explored consensus control under Byzantine attacks, leveraging a virtual twin layer to decouple the defense strategy from the cyber-physical layer and twin layer, which guarantees asymptotic consensus. The Byzantine resilience problem is prevalent in the robotics community \cite{deng2021investigation}, with many works using block-chain technology to secure and sign communications between unaffected robots in a Byzantine-affected multi-robot system \cite{ferrer2021following, strobel2023robot}. With recent advances in learning-based methods, analyzing the effect of Byzantine threats on distributed learning is a critical problem, and has led to statistical machine learning researchers developing Byzantine-robust and -tolerant gradient descent methods \cite{blanchard2017machine, chen2017distributed}.

There are very few existing works that explore enhancing the ADMM algorithm to be suitable for Byzantine- or adversary-resilient problems. The proposed robust ADMM algorithm in \cite{li2022robust} maintains constant communication links between networked nodes, but replaces received information from neighboring nodes if its ``deviation statistics" exceed a threshold, determined by graph properties. The authors of \cite{he2025straggler} explored the straggler-resilience problem by enhancing conventional ADMM to be capable of handling asynchronous nodes. The Byzantine threat problem is complex due to the freedoms that Byzantine nodes have; dealing with adversarial Byzantine attacks is more challenging. As demonstrated in previous works, it is impossible to create a single algorithm that can always compute the optimal solution in the attack-free case while also being resilient in the attack case \cite{sundaram2018distributed}. As such, researchers seek to develop resilient distributed optimization algorithms that produce a solution that is bounded within a neighborhood of the global minima. This provides guarantees on the performance of these algorithms in the face of even adversarially designed attacks. 

In this paper, we develop a resilient ADMM framework that leverages dynamic edge weights to reduce the effect of affected nodes on the distributed optimization process. We guarantee the existence of an update rule for the edge weights that ensure the proposed algorithm converges to the global optimal solution in the Byzantine-free scenario and demonstrate the validity of this result through numerical simulations on an example network. Additionally, we derive sufficient conditions for the existence of an update rule for edge weights that results in the proposed ADMM method converging to a bounded stationary point for a network affected by Byzantine nodes.

The rest of the paper is organized as follows. 
Section \ref{sec:formulation} presents the problem formulation, while Section \ref{sec:prop_alg} describes the proposed algorithm.
In Section \ref{sec:bounded_alg}, we develop the main theoretical results of the paper.
Section \ref{sec:sim} presents the numerical simulation results and analysis that support the proposed algorithm and its theoretical findings.
Finally, Section \ref{sec:conclusion} concludes the paper and presents further research directions.

%%% Problem Formulation
\section{Problem Formulation}
\label{sec:formulation}
%   %   %   %   %   %   %   %   %   %   %   %   %
%   NOTATIONS
%   %   %   %   %   %   %   %   %   %   %   %   %
\subsection{Notations}

The set of real numbers and integers are denoted by $\mathbb{R}$ and $\mathbb{Z}$, respectively, and the superscript $+$ stands for the non-negativeness of the set (e.g., $\mathbb{Z}^+$). Then, $\mathbb{R}^n$ is a real column vector of length $n$, and $\mathbb{R}^{n \times m}$ is a real matrix with $n$ rows and $m$ columns. $\mathbf{1} \in \mathbb{R}^N$ is a column vector of size $N$ with all $1$ elements. The cardinality of a given index set $\mathcal{S}$ is denoted with $| \mathcal{S} |$, such that $\mathcal{S} = \{1, 2, \cdots, | \mathcal{S} | \}$.
Let $\langle \cdot , \cdot \rangle$ be the Frobenius inner product operator and $\norm{\cdot}_F$ be the Frobenius norm. For a positive semi-definite matrix $M$ and square matrix $A$ of equal dimensions, let $\norm{A}_M = \sqrt{\langle A, MA \rangle}$. Let $\text{null}(M)$ be the nullspace of matrix $M$. For a column vector $v$, let $\text{span}(v)$ indicate the space spanned by $v$.

%   %   %   %   %   %   %   %   %   %   %   %   %
%   NETWORK MODEL
%   %   %   %   %   %   %   %   %   %   %   %   %
\subsection{Network Model}

We consider an undirected dynamically-weighted network graph $\mathcal{G} (\mathcal{V}, \mathcal{E})$ where $\mathcal{V} = \{1, 2, 3, \dots, N\}$ is the set of all vertices (nodes) and $\mathcal{E} \subset \mathcal{V} \times \mathcal{V}$ is the set of the edges (connections) in the graph. Two nodes $i$ and $j$ are considered neighbors if $(i, j) \in \mathcal{E}$, with the set of neighbor nodes of node $i$ being represented by $\mathcal{N}_i$. 
We let each edge $(i, j) \in \mathcal{E}$ in the network have a time-varying weight $a_{ij}^k \in [0, \infty)$, where $k \in \mathbb{Z}^+$ is a time index. The time-varying weighted adjacency matrix $A^k$ associated with the graph $\mathcal{G}$ is then defined as $A^k = \left[ a^k_{ij} \right] \in \mathbb{R}^{N \times N}$. We do not consider self-loops in the network, i.e., $a^k_{ii} = 0$ for $i \in \mathcal{V}$. Note that $a^k_{ij} = 0$ if $(i,j) \notin \mathcal{E}$ and the matrix $A^k$ is symmetric since the graph is undirected. The time-varying weighted degree matrix $D^k = \left[ d^k_{ij} \right] = \mathbb{R}^{N \times N}$ is a diagonal matrix, where $d^k_{ii}=\sum_{j \in \mathcal{N}_i}a^k_{ij}$ and $d_{ij} = 0$ when $i \neq j$.
The time-varying signed Laplacian matrix $L_-^k$ is then defined as $L_-^k = D^k - A^k$, while the time-varying unsigned Laplacian $L_+^k$ is defined as $L_+^k = D^k + A^k$.

%   %   %   %   %   %   %   %   %   %   %   %   %
%   BYZANTINE THREAT MODEL
%   %   %   %   %   %   %   %   %   %   %   %   %
\subsection{Byzantine Threat Model}

We consider the case of a Byzantine threat (i.e., Byzantine fault or attack) in the network, where $b$ Byzantine nodes are capable of sending arbitrarily false or delayed messages to their neighbors \cite{lamport1982byzantine}. 
Let $\mathcal{B} \subset \mathcal{V}$ denote the set of nodes affected by a Byzantine fault or attack in the network. Recall the affected nodes can send different false messages to its neighbors in the Byzantine threat model. Furthermore, the Byzantine nodes have knowledge of the entire network and can collude with one another. We let the remaining ``honest," or unaffected, nodes in the network be denoted by $\mathcal{R} = \mathcal{V} \setminus \mathcal{B}$. 

%   %   %   %   %   %   %   %   %   %   %   %   %
%   DISTRIBUTED OPTIMIZATION
%   %   %   %   %   %   %   %   %   %   %   %   %
\subsection{Distributed Optimization}

We adopt a commonly used distributed minimization problem in \eqref{eq:or_distri_opti}, where the agents collaboratively minimize the sum of the local objective functions.
\begin{equation}
\label{eq:or_distri_opti}
\begin{split}
    \min_{\{ x_i \}} ~ & \hspace{20pt} \sum_{i \in \mathcal{V}} f_i (x_i) \\
    \text{s.t.} ~ & \hspace{20pt} x_i = x_j ~ \forall (i, j) \in \mathcal{E}.
\end{split}
\end{equation}
We define $x_i$ to be the local solution at node $i$. Then, $X$ is the collective solution which is organized as such: 
$$X = \begin{bmatrix} x_1, x_2, \cdots, x_N \end{bmatrix}^\top \in \mathbb{R}^{N \times n}.$$
A collective objective function $f(X) = \sum_{i \in \mathcal{V}} f_i(x_i)$ is defined as a re-expression of the objective function in \eqref{eq:or_distri_opti}.

Suppose $z_i = x_i + e_i$ for $i \in \mathcal{V}$, where $e_i$ is a non-zero error term for $i \in \mathcal{B}$. Let $Z \in \mathbb{R}^{N \times n}$ and $E \in \mathbb{R}^{N \times n}$ be organized similar to $X$ such that $Z = \begin{bmatrix} z_1, z_2, \cdots, z_N \end{bmatrix}^\top$ and $E = \begin{bmatrix} e_1, e_2, \cdots, e_N \end{bmatrix}^\top$. 

We define $X^* = Z^*$ to be the optimal solution of the Byzantine-threat-free optimization problem. This solution is consensual, so all local solutions are equal and optimal to \eqref{eq:or_distri_opti}:
\begin{align*}
    x_1^* = x_2^* = \dots = x_N^*, \quad \quad 
    z_1^* = z_2^* = \dots = z_N^*.
\end{align*}

%%% Proposed Algorithm
\section{Proposed Algorithm Development}
\label{sec:prop_alg}
%   %   %   %   %   %   %   %   %   %   %   %   %
%   DYNAMIC WEIGHTED ADMM UNDER BYZANTINE THREAT
%   %   %   %   %   %   %   %   %   %   %   %   %

%%%     WEIGHTED ADMM - MATRIX

In this section, we propose a dynamically weighted ADMM (DW-ADMM) framework for Byzantine resilience in distributed optimization. We first refer to the proposed weighted ADMM algorithm from \cite{ling2016weighted}, with the following network-wide update law:
\begin{align*}
% \label{eq:or_weighted_admm_matrix}
\begin{split}
    X^{k+1} &= \argmin_X { f(X) + \langle X, \Lambda^k + D X - (D+A)X^k \rangle}, \\
    \Lambda^{k+1} &= \Lambda^k + (D-A) X^{k+1}.
\end{split}
\end{align*}
where $D \in \mathbb{R}^{N \times N}$ is a diagonal matrix, $A \in \mathbb{R}^{N \times N}$ is a symmetric matrix, and the sequence $\Lambda^k$ represents the dual variable of the optimization problem at iteration $k$.
%%%     WEIGHTED ADMM - NODE
The network update law above can be decomposed by node, giving the update law for node $i \in \mathcal{V}$ in \eqref{eq:or_weighted_admm_node}. The weighting variables $d_{ii}$ and $a_{ij}$ are elements of $D$ and $A$, respectively. 
\begin{align}
\label{eq:or_weighted_admm_node}
\begin{split}
     x_i^{k+1} &= \argmin_{x_i} f_i(x_i) + \langle x_i, \lambda_i^k - d_{ii} x_i^k - \sum_{j \in \mathcal{N}_i} a_{ij} x_j^k \rangle \\ &\hspace{1in} + d_{ii} \norm{x_i}^2, \\
     \lambda_i^{k+1} &= \lambda_i^k + d_{ii}x_i^{k+1} - \sum_{j \in \mathcal{N}_i} a_{ij} x_j^{k+1}.
\end{split}
\end{align}

%%%     MODIFIED WEIGHTED ADMM - NODE
Note that the dual variable $\lambda_i$ tracks the total constraint deviation between the local estimate $x_i$ of agent $i$ and all its neighbor's $j\in \mathcal{N}_i$, resulting in the summation term in \eqref{eq:or_weighted_admm_node}. The dual variable $\lambda_i$ can be expressed as the sum of all $\lambda_{ij}$ for $j\in \mathcal{N}_i$, where $\lambda_{ij}$ is defined only for edge $(i,j) \in \mathcal{E}$. This will become important in discussion of implementation of the proposed algorithm.
We introduce time-varying weights on the edges between agents, which can be expressed in terms of the time-varying weighted degree and adjacency matrices, $D^k = \big[ d_{ij}^k \big]$ and $A^k = \big[ a_{ij}^k \big]$, respectively.

Incorporating these dynamic weights as well as the Byzantine threat vector $e_i$ results in the update rule for node $i$:
\begin{align}
\begin{split}
    &\nabla f_i(x^{k+1}_i) + \lambda_i^k + 2 d^k_{ii} x^{k+1}_i = d^k_{ii} z_i^k + \sum_{j \in \mathcal{N}_i} a^k_{ij} z_j^k, \\
    &\lambda_{i}^{k+1} = \lambda_{i}^k + d^k_{ii} z^{k+1}_i - \sum_{j \in \mathcal{N}_i}a^k_{ij} z^{k+1}_j.
\end{split}
\end{align}

Expressing this update for the whole network results in the update rule in \eqref{eq:prop_admm_matrix_update}. Recall $L_-^k$ and $L_+^k$ are the time-varying signed and unsigned Laplacian matrices, respectively. This form will be used for the analysis of the algorithm.
\begin{subequations}
\label{eq:prop_admm_matrix_update}
\begin{align}
    \label{eq:prop_x_update}
    &\nabla f(X^{k+1}) + \Lambda^k + 2 D^k X^{k+1} = L_+^k Z^k, \\
    \label{eq:prop_lam_update}
    &\Lambda^{k+1} = \Lambda^k + L_-^k Z^{k+1}.
\end{align}
\end{subequations}

Additionally, we define an update rule for the Laplacian matrices below for a positive definite matrix $M^k \in \mathbb{R}^{N \times N}$: 
\begin{align}
\label{eq:laplacian_update}
    \vishnu{ L_{+/-}^{k+1} = L_{+/-}^k M^k = L_{+/-}^0 M_0^k, } 
\end{align}
where $M_0^k = \prod_{i=0}^k M^i$ and $L_{+/-}^0$ is the original unweighted Laplacian matrix of $\mathcal{G}$. For $k \in \mathbb{Z} \setminus \mathbb{Z}^+$, we let $L^k_{+/-} = L^0_{+/-}$ and $M_0^{k} = I$. It is noted that $\text{null}(M^k) = \text{span}(\mathbf{1})$ so $L_-^{k+1}$ maintains the properties of the signed Laplacian matrix.

Our approach can thus be organized into 3 major steps:
\begin{enumerate}
    \item $X$-update: Solution update $X^{k+1}$ is computed with \eqref{eq:prop_x_update}.
    \item $\Lambda$-update: Dual variable $\Lambda^{k+1}$ is updated using \eqref{eq:prop_lam_update}.
    \item $L$-update: New weighted Laplacian matrices $L_{+/-}^{k+1}$ are found according to \eqref{eq:laplacian_update}.
\end{enumerate}

%%% Boundedness Result
\section{Main Results}
\label{sec:bounded_alg}

%   %   %   %   %   %   %   %   %   %   %   %   %
%   THEORETICAL ANALYSIS
%   %   %   %   %   %   %   %   %   %   %   %   %
\subsection{Problem Assumptions and Supporting Lemmas}

We first provide assumptions about the network $\mathcal{G}$ and the objective function $f(\cdot)$ that are necessary for the analysis of the proposed algorithm. 
\vishnu{Assumption \ref{as:non_bipartite} limits the structure of the graph to be non-bipartite, which provides properties on $L^k_+$, but is not incredibly restrictive in practice as it simply requires the graph to have an odd-length cycle.}
\begin{assumption}
    \label{as:non_bipartite}
    The graph $\mathcal{G}$ is non-bipartite, i.e., the nodes cannot be partitioned into two sets such that no two nodes from the same set are adjacent \cite{asratian1998bipartite}. 
\end{assumption}

\vishnu{Assumption \ref{as:honest_nodes} requires the non-Byzantine nodes to be connected and synchronized. Similar assumptions have been made in previous works on Byzantine resilience \cite{su2020byzantine}.}

\begin{assumption}
    \label{as:honest_nodes}
    The honest nodes in the network $\mathcal{R} \subset \mathcal{V}$ form a connected subgraph $\mathcal{H}$ of $\mathcal{G}(\mathcal{V}, \mathcal{E})$. Additionally, the actions of the honest nodes in the network $\mathcal{R}$ are synchronized.
\end{assumption}

\vishnu{Assumptions \ref{as:solution_set} and \ref{as:convex_cont_diff} require the objective function to have a solution, be convex, and be continuously differentiable. These are common assumptions in distributed optimization \cite{li2022robust, ling2016weighted, su2020byzantine, sundaram2018distributed}.}

\begin{assumption}
    \label{as:solution_set}
    The solution set to \eqref{eq:or_distri_opti} is nonempty and contains at least one bounded element, denoted using $\mathcal{X}^*$.
\end{assumption}

\begin{assumption}
    \label{as:convex_cont_diff}
    The local objective functions $f_i(\cdot)$ are convex and continuously differentiable.
\end{assumption}

%   %   %   %   %   %   %   %   %   %   %   %   %
%   USEFUL THINGS FOR LATER
%   %   %   %   %   %   %   %   %   %   %   %   %
We now present lemmas that are vital to the theoretical analysis of the proposed algorithm. First, we define a new series of matrices $Y^k \in \mathbb{R}^{N \times n}$
related to the dual variable $\Lambda^k$ and show that its update rule is equivalent to that in \eqref{eq:prop_lam_update}.
\begin{lemma} \label{lem:new_dual_var_update}
    The following are equivalent to the dual variable update step in \eqref{eq:prop_lam_update}, where $N = \big( L_-^0 \big)^{1/2}$ and $Y^0 = 0$:
    \begin{subequations}
    \begin{align}
        \label{eq:l1_2}
        Y^{k+1} &= Y^k + N M_0^{k-1} Z^{k+1}, \\
        \label{eq:l1_1}
        Y^k &= \sum_{i=1}^k N M_0^{i-2} Z^i,
    \end{align}
    \end{subequations}
\end{lemma}
\begin{proof}
    Applying the Laplacian update \eqref{eq:laplacian_update} to \eqref{eq:prop_lam_update} gives:
    \begin{equation*}
        \Lambda^{k+1} = \Lambda^k + L_-^0 M_0^{k-1} Z^{k+1}.
    \end{equation*}
    It is clear this equation is equivalent to \eqref{eq:l1_2} by the identity $\Lambda^k = N Y^k$. 
    $Y^k$ is guaranteed to exist
    for a given $\Lambda^k$ for all $k \in \mathbb{Z}^+$ when $\Lambda^k$ is in the column space of $L_-^0$ and $N$. Since $\Lambda^0$ is initialized as $0$, it is in the column space of $L_-^0$ and $N$, and all subsequent $\Lambda^k$ remain in their column space due to the recursion shown above. 
    Evaluating the recursion as a summation results in \eqref{eq:l1_1}.
\end{proof}

With a new update rule for $\Lambda^k$ written, we now rewrite the update step for $Z^{k+1}$.

\begin{lemma} \label{lem:Zk_sequence}
    The sequences of optimization variables must adhere to the following:
    \begin{align}
    \begin{split}
        \label{eq:l2}
        L_+^k ( Z^{k+1} - Z^k ) - 2 D^k E^{k+1} = -\Lambda^{k+1} - \nabla f(X^{k+1}).
    \end{split}
    \end{align}
\end{lemma}
\begin{proof}
    Using \eqref{eq:prop_admm_matrix_update}, we can write the following relation:
    \begin{align}
        2 D^k X^{k+1} - L_+^k Z^k &= - \Lambda^{k} - \nabla f(X^{k+1}).
    \end{align}
    Subtracting $L_-^k Z^{k+1}$ from both sides and rewriting the left expression: 
    \begin{align}
        2 D^k X^{k+1} - L_+^k Z^k - L_-^k Z^{k+1}
        &= - \Lambda^{k+1} - \nabla f(X^{k+1}), \nonumber\\
        = L_+^k (Z^{k+1}& - Z^k) - 2 D^k E^{k+1}.
    \end{align}
    Then, we have the desired result.
\end{proof}

We now define the first-order optimality condition for the problem. This will indicate where the optimization sequences have reached an optimal point. 
\begin{lemma} \label{lem:first_order_opt}
    The optimal solution $(X^*, Y^*)$ satisfies the following for $N = {\big( L_-^0 \big)}^{1/2}$ and $X^* = Z^*$:
    \begin{subequations}
    \label{eq:l3}
    \begin{align}
        \label{eq:l3_1}
        0 &= \nabla f(X^*) + N Y^*, \\
        \label{eq:l3_2}
        0 &= L_-^0 X^* = L_-^k X^*,
    \end{align}
    \end{subequations}
\end{lemma}
\begin{proof}
    The first-order optimality condition is derived with respect to the optimal solution $X^* = Z^*$ from conventional ADMM, where the Laplacian matrices are constant. That is,  $L^k_{+/-} =L_{+/-}^0$ for $\forall k \in \mathbb{Z}^+$.
    Since $X^* = Z^*$ is consensual and $\text{null}(L_-^0) = \text{span}(\mathbf{1})$ by the property of the signed Laplacian matrix, the second equation is derived. 

    Since $X^*$ is a consensual optimal solution, we can also write $\mathbf{1}^\top \nabla f(X^*) = 0$, which is equivalent to saying $\nabla f(X^*)$ stays in the span of $L_-^0$. That is, $\nabla f(X^*) + NY^* = 0$, where $Y^* = NP$ for some $P\in \mathbb{R}^{N \times n}$.
\end{proof}

Finally, based on the earlier lemmas, we define the recursion of the proposed algorithm.
\begin{lemma} Recursion of the proposed ADMM algorithm:
    \begin{align}
    \begin{split}
        \label{eq:l4}
        \nabla f(X^{k+1}) &= \nabla f(X^*) -( \Lambda^{k+1} - \Lambda^* ) \\ & \hspace{20pt} + L_+^k(Z^k - Z^{k+1}) + 2 D^k E^{k+1}, \\
        Y^{k+1} &= Y^k + N M_0^{k-1} (Z^{k+1} - Z^*).
    \end{split}
    \end{align}
\end{lemma}
\begin{proof}
    We leverage Lemmas \ref{lem:new_dual_var_update}, \ref{lem:Zk_sequence}, and \ref{lem:first_order_opt} to derive \eqref{eq:l4}. Subtracting \eqref{eq:l3_2} from \eqref{eq:l1_2}, the dual variable update is: 
    \begin{equation}
        Y^{k+1} = Y^k + N M_0^{k-1} Z^{k+1} - L_-^0 Z^*.
    \end{equation}
    Since $L_-^k$ retains the properties of the signed Laplacian matrix, i.e., $\text{null}(L_-^k) = \text{span}(\mathbf{1})$,
    it is true that $L_-^k Z^* = L_-^0 Z^* = 0$. Applying this identity with the Laplacian update rule \eqref{eq:laplacian_update} gives the desired result for $Y^{k+1}$.
    By \eqref{eq:l2} and \eqref{eq:l3_1}, it yields that
    \begin{align}
        \nabla f(X^{k+1}) &= \nabla f (X^*) - (\Lambda^{k+1} - N Y^*) \nonumber \\ & \hspace{20pt} + L_+^k (Z^k - Z^{k+1}) + 2D^k E^{k+1}.
    \end{align}
    Applying the identity $\Lambda^* = NY^*$ gives the desired result.% for $\nabla f( X^{k+1})$.
\end{proof}

\subsection{Convergence and Boundedness}

We first consider the case where there is no error in the network. We show that there exists a sequence of Laplacian update matrices $M^k$ such that the proposed algorithm will converge to the global optimal solution.

\begin{theorem}[Error-Free Convergence]
\label{thm:1}
    Consider the distributed optimization problem in \eqref{eq:or_distri_opti} and suppose that the objective function and the underlying network of the nodes satisfy Assumptions \ref{as:non_bipartite} - \ref{as:convex_cont_diff}. In the error-free case, where $E^k = 0$ $\forall k \in \mathbb{Z}^+$, with $M^k$ being a function of the iterate variables at iteration $k$, $\alpha^k I \succeq M^k \succ 0$ for $\alpha^k \in \mathbb{R}$, and $\text{null}(M^k) = \text{span}(\mathbf{1})$ $\forall k \in \mathbb{Z}^+$, the algorithm converges to the global optimal solution $X^*$.
\end{theorem}
\begin{proof}
    We first use the convexity property of $f(X)$ supposed in Assumption \ref{as:convex_cont_diff} before setting $E^k = 0$:
    \begin{align*}
        0 &\leq \langle X^{k+1} - X^*, \nabla f(X^{k+1}) - \nabla f(X^*) \rangle,
    \end{align*}
    \vspace{-18pt}
    \begin{align}
    \label{eq:convexity_result}
        \begin{split}
            0 &\leq \langle Z^{k+1} - Z^*, L_+^k(Z^k - Z^{k+1}) \rangle \\ & \quad + \langle Z^{k+1} - Z^*, -(\Lambda^{k+1} - \Lambda^*) \rangle + \mathcal{E}^{k+1},
        \end{split}
    \end{align}
    where $\mathcal{E}^{k+1} = \langle E^{k+1}, L_+^k (Z^{k+1} - Z^k) + (\Lambda^{k+1} - \Lambda^*) + 2D^k (X^{k+1} - X^*)\rangle$.
    Now setting $E^k = 0$ for $\forall k \in \mathbb{Z}^+$:
    \begin{align*}
        0 &\leq \langle X^{k+1} - X^*, L_+^k(X^k - X^{k+1}) \rangle \\ 
        & \quad + \langle X^{k+1} - X^*, -(\Lambda^{k+1} - \Lambda^*) \rangle.
    \end{align*}
    Rewriting the expression on the right hand side:
    \begin{align*}
        \langle X^{k+1} - X^*, L_+^k(X^k - X^{k+1}) \rangle + \langle X^{k+1} - X^*, \\ -(\Lambda^{k+1} - \Lambda^*) \rangle
        = \langle Q^{k+1} - Q^*, G^k (Q^k - Q^{k+1}) \rangle,
    \end{align*}
    where
    \begin{equation*}
        Q^k = \begin{bsmallmatrix} Y^k \\ X^k \end{bsmallmatrix}, \quad 
        G^k = \begin{bsmallmatrix} (M_0^{k-1})^{-1} & 0 \\ 0 & L_+^k \end{bsmallmatrix}.
    \end{equation*}

    This can be written as:
    \begin{align*}
        \langle Q^{k+1} - Q^*, G^k (Q^k - Q^{k+1}) \rangle = \frac{1}{2} \Big[ \norm{Q^k - Q^*}^2_{G^k} \\ - \norm{Q^{k+1} - Q^*}^2_{G^k} - \norm{Q^k - Q^{k+1}}^2_{G^k} \Big]
        & \geq 0.
    \end{align*}
    Summing the above over $k \in \mathbb{Z}^+$ results in a telescoping sum:
    \begin{align*}
        0 \leq \norm{Q^0 - Q^*}_{G^0}^2 - \norm{Q^\infty - Q^*}_{G^\infty}^2 - \sum_{k=0}^\infty {\norm{Q^k - Q^{k+1}}_{G^k}^2}, \\
        \sum_{k=0}^\infty {\norm{Q^k - Q^{k+1}}_{G^k}^2} \leq \norm{Q^0 - Q^*}_{G^0}^2 - \norm{Q^\infty - Q^*}_{G^\infty}^2.
    \end{align*}
    Since both $\norm{Q^0 - Q^*}_{G^0}^2$ and $\norm{Q^\infty - Q^*}_{G^\infty}^2$ are bounded, the summation is also bounded and
    \begin{align*}
        \lim_{k \rightarrow \infty} \norm{Q^k - Q^{k+1}}_{G^k}^2 = 0,
    \end{align*}
    which, by the definition of $Q^k$ and $G^k$, implies
    \begin{align*}
        \lim_{k \rightarrow \infty} \norm{Y^k - Y^{k+1}}_{(M_0^{k-1})^{-1}}^2 = 0, \hspace{4pt} \lim_{k \rightarrow \infty} \norm{X^k - X^{k+1}}_{L_+^k}^2 = 0.
    \end{align*}

    We first evaluate the limit with the $Y$ terms. Since $M^k$ is spectrally upper-bounded, i.e., $\alpha^k I \succeq M^k$, it is also true that $\alpha I \succeq M_0^{k-1}$, where $\alpha = \prod_{0}^{k-1} \alpha^k$, and $\frac{1}{\alpha}I \preceq (M_0^{k-1})^{-1}$. For a matrix $M$ and symmetric matrix $G$, $\lambda_\text{min} (G) \norm{M}_F^2 \leq \norm{M}_G^2$ \cite{boyd2004convex}. Then, we can show that
    \begin{align*}
        {\alpha}^{-1} \norm{Y^k - Y^{k+1}}_F^2 \leq \norm{Y^k - Y^{k+1}}_{(M_0^{k-1})^{-1}}^2.
    \end{align*}
    Applying the limit as $k \rightarrow \infty$ gives
    \begin{align*}
        {\alpha}^{-1} \lim_{k \rightarrow \infty} \norm{Y^k - Y^{k+1}}_F^2 &\leq \lim_{k \rightarrow \infty}  \norm{Y^k - Y^{k+1}}_{(M_0^{k-1})^{-1}}^2, \\
        \lim_{k \rightarrow \infty} \norm{Y^k - Y^{k+1}}_F^2 &= 0.
    \end{align*}
    This shows that $Y^k$ converges to a stationary point $\hat{Y}$. We now evaluate the limit with the $X$ terms. From \eqref{eq:laplacian_update}, it can be seen that $L_+^k=L_+^0 M_0^{k-1}$, which is positive definite by Assumption \ref{as:non_bipartite}. The positive definiteness means that $L_+^k$ is spectrally lower-bounded, $L_+^k \succeq \beta I$ for $\beta \in \mathbb{R}$, and 
    \begin{align*}
        \beta \norm{X^k - X^{k+1}}_F^2 &\leq \norm{X^k - X^{k+1}}_{L_+^k}^2, \\
        \beta \lim_{k \rightarrow \infty} \norm{X^k - X^{k+1}}_F^2 &\leq \lim_{k \rightarrow \infty} \norm{X^k - X^{k+1}}_{L_+^k}^2, \\
        \lim_{k \rightarrow \infty} \norm{X^k - X^{k+1}}_F^2 &= 0,
    \end{align*}
    showing $X^k$ converges to a stationary point $\hat{X}$.
    We substitute the stationary pair $(\hat{X}, \hat{Y})$ into the recursion formulas
    \eqref{eq:prop_x_update} and \eqref{eq:l1_2}. Recall that $M^k$, and thus $M_0^{k-1}$, are functions of the optimization iterate variables. Since the optimization has reached a stationary point, we can conclude that $M_0^{k-1}$ has also reached a stationary point $\hat{M}$, leading to stationary degree and adjacency matrices, $\hat{D}$ and $\hat{A}$, respectively. This gives
    \begin{align*}
        & \nabla f(\hat{X}) + N \hat{Y} + \hat{L}_- \hat{X} = 0, \\
        & 0 = N \hat{M} \hat{X}.
    \end{align*}
    This is equivalent to the optimality condition \eqref{eq:l3}, 
    thus showing that the stationary pair is optimal $(X^*, Y^*)$.
\end{proof}

Using Theorem \ref{thm:1}, we have demonstrated that there exists a $\{M^k\}_{k=0}^\infty$ for which our proposed resilient ADMM framework converges to the optimal solution in the absence of an error caused by Byzantine nodes in the network. In Section \ref{sec:sim}, we explore an update rule based on the optimization iterate variables and demonstrate the algorithm's convergence to the optimal solution. 
\vishnu{It must be noted that network topology plays a role in the algorithm's performance, notably the non-bipartite nature of the graph network, which can be preserved by ensuring that a cycle of odd-length exists in the network. Conclusions on transient behavior, such as convergence rates, may be determined with further assumptions on the problem, such as strong convexity of the objective function in \eqref{eq:or_distri_opti}.}

We now consider the case where there are Byzantine nodes affecting the network and show that the proposed method can produce a bounded result. \vishnu{We propose integrating our DW-ADMM framework with trust computation methods such as those described in \cite{govindan2011trust, zegers2021event}, where each node calculates a \textit{trust level} of its neighbors, typically defined from $0$ to $1$, which represents the level of confidence that a neighboring node is honest. These trust levels can be computed by measuring discrepancies in received information over each network edge. In this work, when trust in a node drops below a defined level, the node is identified as a faulty/attacked node. 
Location and sparsity of the Byzantine nodes may affect convergence rates.}

\begin{theorem}[Boundedness under Byzantine Threat]
\label{thm:2}
    Consider the distributed optimization problem in \eqref{eq:or_distri_opti} and suppose that the objective function and the underlying network of the agents satisfy Assumptions \ref{as:non_bipartite} - \ref{as:convex_cont_diff}. 
    Let $M^k$ be a function of the iterate variables at iteration $k$, $\alpha^k I \succeq M^k \succ 0$ for $\alpha^k \in \mathbb{R}$, and $\text{null}(M^k) = \text{span}(\mathbf{1})$ $\forall k \in \mathbb{Z}^+$.
    Suppose Byzantine nodes $j \in \mathcal{B} \subset \mathcal{V}$ are detected and identified at iteration $k' \in \mathbb{Z}^+$. An $M^{k'}$ can be designed such that weight of the edges incident to the nodes become $0$, resulting in the unaffected nodes converging to a stationary pair $(\hat{X}, \hat{Y})$. 
\end{theorem}
\begin{proof}
    We start with \eqref{eq:convexity_result} from Theorem \ref{thm:1}, but $E^k$ will not be set to $0$ as we consider the error-present case. Similar methods from Theorem \ref{thm:1} will be used to simplify expressions.
    \begin{align*}
        \begin{split}
            &0 \leq \mathcal{E}^{k+1} + \langle Z^{k+1} - Z^*, L_+^k(Z^k - Z^{k+1}) \rangle + \langle Z^{k+1} - Z^*, \\ 
            & \quad -(\Lambda^{k+1} - \Lambda^*) \rangle
            \leq \frac{1}{2} \Big[ \norm{Q^k - Q^*}^2_{G^k} - \norm{Q^{k+1} - Q^*}^2_{G^k} \\
            & \quad \quad - \norm{Q^k - Q^{k+1}}^2_{G^k} \Big] + \mathcal{E}^{k+1}
        \end{split}
    \end{align*}
    where
    \begin{equation*}
        Q^k = \begin{bsmallmatrix} Y^k \\ Z^k \end{bsmallmatrix}, \quad 
        G^k = \begin{bsmallmatrix} (M_0^{k-1})^{-1} & 0 \\ 0 & L_+^k \end{bsmallmatrix}.
    \end{equation*}

    Summing the above over $k \in \mathbb{Z}^+$:
    \begin{align*}
        \sum_{k=0}^\infty {\norm{Q^k - Q^{k+1}}_{G^k}^2} &\leq \norm{Q^0 - Q^*}_{G^0}^2 - \norm{Q^\infty - Q^*}_{G^\infty}^2 \\
          & \quad \quad + 2 \sum_{k=0}^\infty \mathcal{E}^{k+1} .
    \end{align*}
    
    Recall that a Byzantine threat is detected at $k'$ and $M^{k'}$ is designed such that the edge weights associated with the identified Byzantine nodes become $0$. Then, the summation of $\mathcal{E}^{k+1}$ becomes: $\sum_{k=0}^\infty \mathcal{E}^{k+1} = \sum_{k=0}^{k'} \mathcal{E}^{k+1}$,
    which is a bounded expression. Then, 
    since both $\norm{Q^0 - Q^*}_{G^0}^2$ and $\norm{Q^\infty - Q^*}_{G^\infty}^2$ are bounded, the summation $\sum_{k=0}^\infty {\norm{Q^k - Q^{k+1}}_{G^k}^2}$ is also bounded and $\lim_{k \rightarrow \infty} \norm{Q^k - Q^{k+1}}^2_{G^k} = 0$.
    By definition of $Q^k$ and $G^k$, this implies
    \begin{align*}
        \lim_{k \rightarrow \infty} \norm{Y^k - Y^{k+1}}_{(M_0^{k-1})^{-1}}^2 = 0, \hspace{6pt} \lim_{k \rightarrow \infty} \norm{Z^k - Z^{k+1}}_{L_+^k}^2 = 0.
    \end{align*}

    Note that these limits can be rewritten as below
    \begin{align*}
        & \lim_{k \rightarrow \infty} \norm{ \big( M_0^{k-1} \big)^{-1/2} \big(Y^k - Y^{k+1} \big) }_F^2 = 0, \\
        & \lim_{k \rightarrow \infty} \norm{ \big( {L_+^0 M_0^{k-1}} \big)^{1/2} \big( Z^k - Z^{k+1} \big) }_F^2 = 0.
    \end{align*}
    It can be seen that all error terms are multiplied by a form of $M_0^{k-1}$, and thus $M^k$, so a properly designed update law can eliminate the error terms $E^k$ and $E^{k+1}$. This readily implies the difference between the iterate variables of the unaffected nodes approach zero, and therefore a stationary point.
\end{proof}

%%% Numerical Simulation
\section{Simulation Results}
\label{sec:sim}
\begin{figure}[t] 
\centering
  \includegraphics[width=0.33\textwidth, keepaspectratio]{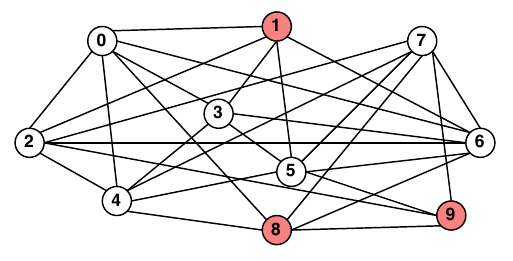}
  \caption{Network Topology for Simulations}
  \label{fig:nettop}
\end{figure}

We demonstrate the performance of the proposed algorithm in numerical simulations with $10$ nodes ($|\mathcal{V}| = 10$), comparing the results with conventional ADMM. The network topology is shown in Figure \ref{fig:nettop}, with the chosen Byzantine nodes indicated in red. First, we show that the proposed resilient ADMM and conventional ADMM produce nearly identical results under certain choices for the Laplacian update in the error-free case, validating Theorem \ref{thm:1}. 
We then show that the proposed algorithm produces a bounded solution when affected by Byzantine nodes, supporting Theorem \ref{thm:2}. The update rules for the time-varying Laplacian matrices are designed to be functions of the optimization iterate, dual variable, and an auxiliary variable $\Phi^k$ that sums the total unsigned constraint deviation, i.e., $\Phi^{k+1} = \Phi^k + L^k_-  \big( Z^k \big) ^{\circ 2}$, where $ \big( Z^k \big)^{\circ 2}$ denotes the element-wise square of the matrix $Z^k$. We define a scaling parameter $\alpha \in \mathbb{R}$ and a weight threshold $\tau \in \mathbb{R}$, and let the edge weights $a_{ij}^k$ be updated as follows
\begin{figure*}[t] 
\centering
  \includegraphics[width=0.33\textwidth, keepaspectratio]{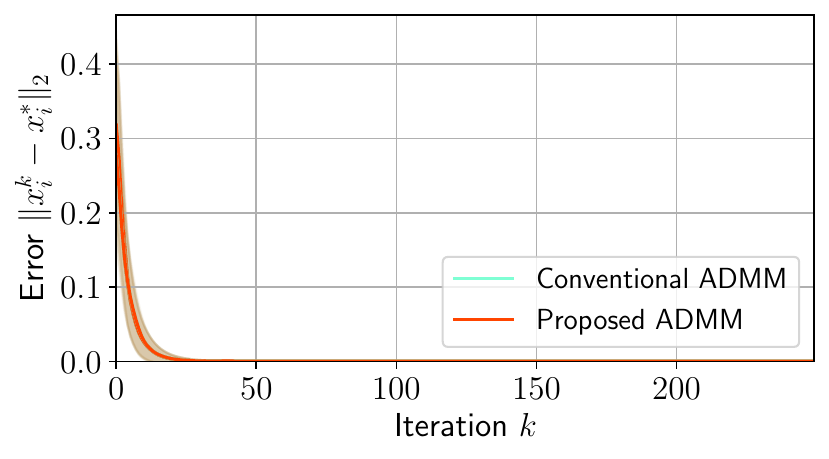}
  \includegraphics[width=0.32\textwidth, keepaspectratio]{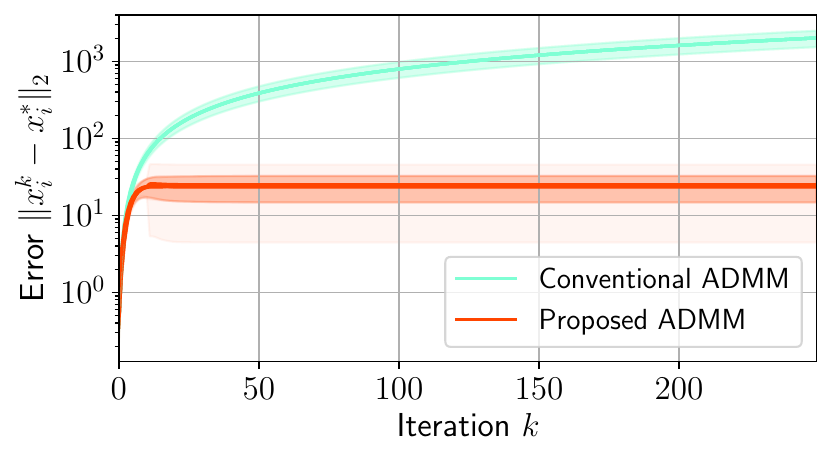}
  \includegraphics[width=0.33\textwidth, keepaspectratio]{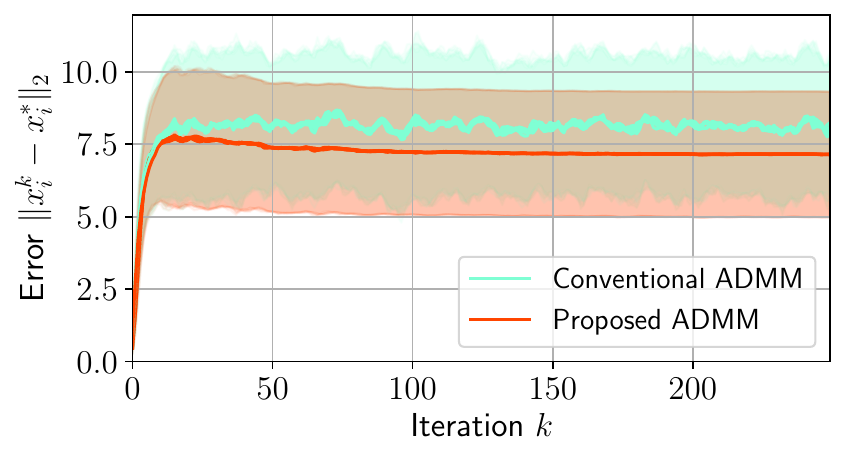}
  \caption{\textbf{Left} Convergence without Error, \textbf{Center} Boundedness with Ramp Error, \textbf{Right} Boundedness with Random Error.
  }
  \label{fig:mc_conv}
\end{figure*}
\begin{align}
    \begin{split}
        a^{k}_{ij} &= 
            \begin{cases} 
                g(\alpha, \tau, \phi^k_{ij}) & \text{if} \hspace{6pt} g(\alpha, \tau, \phi^k_{ij}) > \tau, \\
                0 & \text{otherwise}.
            \end{cases}
    \end{split}
\end{align}
where $g(\alpha, \tau, \phi^k_{ij}) = \frac{\alpha}{\norm{\phi^k_{ij}} + \alpha}$
and $\phi^k_{ij}$ is an element of $\Phi^k$. The edge weight $a^k_{ij} \in [0, 1]$ can be thought of as the trust level between nodes $i$ and $j$ that adaptively changes based on the amount that the consensus constraint deviates. In all simulations, $\alpha = 100$ and $\tau = 0.2$.
For both cases, we consider a quadratic objective function for node $i \in \mathcal{V}$, i.e., $f_i(x_i) = \norm{y_i - N_i x_i}_2^2$,
with $x_i \in \mathbb{R}^{3}$, $y_i \in \mathbb{R}^{3}$, and $N_i \in \mathbb{R}^{3 \times 3}$. The elements of $y_i$ and $N_i$ are generated according to the standard normal distribution $\mathcal{N}(0,1)$. For each simulation described below, both conventional ADMM and the proposed ADMM algorithms are run for $250$ iterations with the same initial guess, with the results averaged over $100$ Monte Carlo (MC) simulations. The results are plotted with the $1$ standard deviation bounds.

We first compare the algorithms in the error-free case. The optimization solution error $\norm{x_i^k - x_i^*}_2$ at each iteration $k$ is plotted in the left side of Fig. \ref{fig:mc_conv}, where it can be seen that the final error of the conventional ADMM and proposed ADMM are approximately equal, with a difference on the order of $10^{-4}$. Furthermore, both algorithms converge quickly to the minimizer, indicating that there is not any marked difference between the two algorithms supporting our theoretical results.

Next, we test the proposed algorithm under two error-present cases: (i) the Byzantine nodes add a linear ``ramp" error function (e.g., $e_i^k = 1 + k$) to their updates and (ii) the Byzantine nodes add Gaussian noise (e.g., $e_i^k \in \mathcal{N}(1, 1.5^2)$) to their optimization updates. The results from simulating case (i) can be seen in the center of Fig. \ref{fig:mc_conv}, which shows the proposed algorithm converging to a stationary point, while the solution from conventional ADMM continues to grow with $e_i^k$. Note that the y-axis is plotted with a log-scale, so the difference between the final solutions are on the order of $10^{3}$. The right of Fig. \ref{fig:mc_conv} shows the results from case (ii). The difference between the final errors of the proposed and conventional methods is not as significant as that in case (ii). However, it is clear that the solution from proposed ADMM algorithm is closer to the true solution than that from the conventional ADMM algorithm, while also converging to a stationary point.

%%% Conclusion
\section{Conclusion}
\label{sec:conclusion}

We have designed the DW-ADMM framework that uses dynamic edge weighting to promote resilience in the presence of Byzantine nodes in a network. We showed that the proposed framework can be designed to converge to the true optimal solution in an error-free scenario, and also converge to a bounded stationary point in the presence of Byzantine nodes in the network. These results were demonstrated and supported with rigorous Monte Carlo simulations.
In future work, we plan to \vishnu{investigate the convergence rate of DW-ADMM under strong convexity assumptions}. We also plan to explore implementing the proposed resilient ADMM algorithm on a multi-robot system, and devise a hierarchical mechanism to utilize the edge weights to restructure the system topology while maintaining a resilient structure.

%%% Bibliography
\bibliographystyle{ieeetr}
\bibliography{IEEEabrv,bibs} 

\end{document}